\documentclass[11pt]{article}
\usepackage{mathrsfs}
\usepackage{latexsym,amsmath,amssymb,amsfonts,epsfig,cite,psfrag,hyperref}
\usepackage{eepic,color,colordvi,amscd,indentfirst,graphics}
\usepackage[ruled,linesnumbered]{algorithm2e}
\usepackage{epsfig}
\usepackage{epstopdf}

\newtheorem{theo}{Theorem}[section]

\newtheorem{lem}[theo]{Lemma}
\newtheorem{coro}[theo]{Corollary}

\def\qed{\hfill \rule{4pt}{7pt}}
\def\pf{\noindent {\it Proof.} }

\begin{document}

\title{The point-thickness of complete multipartite graphs }
\author{Wenzhong Liu\footnote {Corresponding author, E-mail:
wzhliu7502@nuaa.edu.cn.}, \ Wangkai Zhang\footnote{E-mail:
zhangwangkai257@163.com.} \\
School of Mathematics,
Nanjing University of Aeronautics and Astronautics,\\
Nanjing 210016, P.R. China}

\date{}

\maketitle{}

\begin{abstract}
 The point-thickness $\theta'(G)$ of a graph $G$ is the minimum number of subsets into which the vertex set $V(G)$ of $G$
 is partitioned such that each subset induces a planar subgraph. In this paper, we determine the  point-thickness of
 complete multipartite graphs. As a special case, we also obtain the  point-thickness of complete graphs.
\medskip

\smallskip


\end{abstract}

\section{Introduction}

A graph $G$ is a pair $(V(G), E(G))$ with vertex set $V(G)$ and edge set $E(G)$.
The thickness  $\theta(G)$ of a graph $G$ is the minimum number of planar subgraphs into
which G can be decomposed.
Since Tutte \cite{Tutte} introduced the thickness of a graph,  some classes of  graphs whose
thickness was determined are complete graphs $K_{n}$\cite{AG, BH},  almost all complete bipartite graphs $K_{m, n}$\cite{BHM},
 hypercube graphs $Q_n$\cite{MK}. Determining 
 the thickness of a graph is $NP$-hard \cite{MA}. For many graphs,  attention has been focused on  upper bounds of  their thickness, see for example \cite{AS,CY, DH, MD, JMOS, XZ}. One can refer to  \cite{MOS} for more  on the thickness of graphs.

As a vertex version of the thickness, the point-thickness $\theta'(G)$ of a graph $G$ is the minimum number of subsets into which the vertex set $V(G)$
 is partitioned such that each subset induces a planar subgraph.  Considering induced subgraphs of a vertex partition of a graph is a classical problem in graph theory.  When
replacing a planar subgraph in the point-thickness definition by an independent set, a forest, and a k-degenerate graph,
 the minimum number is the chromatic number $\chi(G)$, the point-arboricity $\rho(G)$
 and the point partition number  $\rho_k(G)$ of a graph, respectively.   It is easily checked that
 $\chi(G)\geq \rho(G) \geq \theta'(G)\geq \rho_5(G)$ according to their definitions.
 For the chromatic number $\chi(G)$,  many nice results have been obtained.
 Chartrand et al \cite{CK, CKW} provided formulas for the point-arboricity of complete graphs and complete multipartite graphs, and proved that the point-arboricity of a planar graph is less than 3. Lick and White \cite{LW}  systematically investigated the point partition number of a graph.
There is not any direct result considering the point-thickness of a graph as far as we know.
In fact, each of the four minimum numbers mentioned above is a coloring number, since it gives the minimum number of colorings in any vertex coloring of a graph
 such that each coloring class induces a subgraph with some certain property.
 Related problems of vertex partitions have been studied extensively. 
Borodin \cite{BO} conjectured that every planar graph admits a vertex partition whose parts induce a 1-degenerate and a
2-degenerate graph, respectively, and another vertex partition whose parts induce a 0-degenerate and
a 3-degenerate graph, respectively.  Thomassen \cite{TH95, TH01} proved the two conjectures hold.

In the paper, we focus on the exact value of the point-thickness of graphs,
and determine the  point-thickness of 
 complete multipartite graphs. As a special case, we also obtain the point-thickness of complete graphs.
Before our main theorems are stated, we need to introduce some terminology and notations.
A \emph{$k$-independent-set}, abbreviated to $I_k$-set, is a set of $k$ vertices in which
no two vertices are adjacent.
A \emph{complete $n$-partite graph $K_{p_1, p_2, \cdots, p_n}$} has its vertex set $V$ partitioned
into subsets $V_i$ with $|V_i|=p_i$, $i=1,2, \cdots, n$;  two vertices $x$ and $y$ are adjacent if and only if
$x \in V_i$ and $y\in V_j$ for $i\neq j$.  Without loss of generality, we always assume  $p_1\leq  p_2\leq \cdots \leq p_n$
(if necessary, we can change the orders of some $p_i$ ).
If there are $k$ numbers $ p_{i+1},  p_{i+2}, \cdots,  p_{i+k}$ satisfying
$p_{i+1}=\cdots= p_{i+k}=s$, then $K_{p_1, p_2, \cdots, p_n}$ is briefly denoted by
$K_{p_1, p_2, \cdots, p_{i}, s^{k}, p_{i+k+1}, \cdots  p_{n}}$.
In particular,  if $p_1=p_2=\cdots=p_n=1$,  then
 $K_{p_1, p_2, \cdots, p_n}$ is a \emph{complete graph} on $n$ vertices, also denoted by $K_n$.

   Let $k_1$ and $k_2$ be any two nonnegative integers, we define

   $$ N(k_1, k_2)=\left\{
                   \begin{array}{ll}
                    k_2+\lceil \frac{1}{4}(k_1-3k_2) \rceil, & \hbox{if $k_1\geq 3k_2$;} \\  \medskip

                   \lfloor   \frac{k_1}{3}\rfloor +\lfloor \frac{1}{3}(k_2- \lfloor  \frac{k_1}{3}\rfloor) \rfloor +\sigma, & \hbox{otherwise,}
                   \end{array}
                 \right.$$
  \medskip

where $$\sigma=\left\{
                \begin{array}{ll}
                  0, & \hbox{if $k_1\equiv 0$(mod 3) and $k_2-\lfloor \frac{k_1}{3}\rfloor\equiv 0 $(mod 3);} \\

                  2, & \hbox{if $k_1\equiv 2$(mod 3) and $k_2-\lfloor \frac{k_1}{3}\rfloor\equiv 2 $(mod 3);} \\

                  1, & \hbox{otherwise.}
                \end{array}
              \right. $$

\begin{theo}\label{theo1.1}
Suppose that $k_1, k_2$, $k_3$ and $n$ are any nonnegative integers
and that $p_0=k_1+2k_2+3k_3$. Let $G$ be a complete multipartite graph $K_{1^{k_{1}}, 2^{k_{2}}, 3^{k_{3}}, p_1, p_2, \cdots, p_n}$
with $4 \leq p_1\leq  p_2\leq \cdots \leq p_n$.  We have

$(a)$~~  if $p_0 \leq 2n$, then
$$\theta'(G)=n-max\{ j| \sum^{j} _{i=0} p_i \leq 2(n-j) \}. $$

$(b)$~~ if $p_0 > 2n$, then

$$\theta'(G)=\left\{
                \begin{array}{ll}
                  n+N(k_1+k_3-2n, k_2+k_3), & \hbox{if $k_1+k_3 \geq 2n$;} \\
                  n+N(\varepsilon, \frac{1}{2}(p_0-2n-\varepsilon)), & \hbox{otherwise,}
                  \end{array}
              \right. $$

where  $k_1+k_3\equiv \varepsilon $ (mod 2).
\end{theo}

From Theorem \ref{theo1.1}, we obtain the following corollary.

\begin{coro}\label{coro1.2}
For a  complete graph $K_n$, then

$$\theta'(K_n)=\lceil   \frac{n}{4}  \rceil.$$

\end{coro}

\section{Main results}

This section starts from the following lemmas.

\begin{lem}\label{lem1}
If $G$ is a complete multipartite graph with $n\ge 7$ vertices other than $K_{1, n-1}$, $K_{2, n-2}$, and $K_{1, 1, n-2}$, then
$G$ is a non-planar graph.
\end{lem}

\pf  It can easily be deduced from the fact that all such  graphs contain a subgraph $K_{3, 3}$.

\bigskip

Let $k_1, k_2$, $k_3$ and $n$ be any four nonnegative integers,
and let $p_0=k_1+2k_2+3k_3$. Suppose that $G$ is a complete multipartite graph $K_{1^{k_{1}},~2^{k_{2}},~3^{k_{3}},~p_1,~p_2, \cdots,~p_n}$
with $4 \leq p_1\leq  p_2\leq \cdots \leq p_n$ and that  the last $n$ subsets $V_1$, $V_2$, $\cdots$, $V_n$ of its vertex partition  satisfy
$|V_i|=p_i$ ($1\le i\le n $).

\begin{lem}\label{lem2}
Let $G, p_0$ and $n$ be defined as above. If $p_0 \leq 2n$, then
$$\theta'(G)=n-max\{ j| \sum^{j} _{i=0} p_i \leq 2(n-j) \}. $$

\end{lem}

\pf We apply induction on $n$. For $n=1$, 
 since $p_0 \leq 2n=2$,  $G$ is $K_{2, p_1}$,  $K_{1, 1, p_1}$, $K_{1, p_1}$ or an $I_{p_1}$-set, 
 which  are all  planar graphs.
 Thus $$\theta'(G)=1=n-max\{ j| \sum^{j} _{i=0} p_i \leq 2(n-j) \}.$$ 

Assume that the formula holds for $n\ge 1$,
and consider the complete multipartite graph $G=K_{1^{k_{1}},~2^{k_{2}},~3^{k_{3}},~p_1,~p_2, \cdots,~p_{n+1}}$.
Let $G_1=K_{1^{k_{1}},~2^{k_{2}},~3^{k_{3}},~p_1,~p_2, \cdots,~p_{n}}$. 
Clearly, $G_1$ is a subgraph of $G$. For the graph $G_1$,
suppose that
$$\sum^{t} _{i=0} p_i \leq 2(n-t)~~~ and ~~~ \sum^{t+1} _{i=0} p_i > 2(n-t-1). ~~~~~~~~~~~~(*) $$
 By hypothesis $\theta'(G_1)=n-t$.  Since $G_1$ is a subgraph of  $G$, $\theta'(G) \ge \theta'(G_1)=n-t$.
Since the additional set of $p_{n+1}$ vertices used in forming $G$ induces an 
$I_{p_{n+1}}$-set (a planar graph),
 $\theta'(G) \leq \theta'(G_1)+1$.
We now discuss two cases below.

\medskip
\textit{  Case 1: Suppose $\sum^{t+1} _{i=0} p_i > 2(n+1-(t+1))=2(n-t)$. }
\medskip

By the first inequation in $(*)$, 
$$\sum^{t} _{i=0} p_i \leq 2(n-t)\leq 2(n+1-t),$$  which implies
$$ n+1-max\{ j| \sum^{j} _{i=0} p_i \leq 2(n+1-j) \}=n+1-t=\theta'(G_1)+1.$$
In this case, we need to prove $\theta'(G)=\theta'(G_1)+1$. If not, then $\theta'(G)=\theta'(G_1)=n-t$.
Let $G_2$ be a complete multipartite graph $K_{1^{k_{1}},~2^{k_{2}},~3^{k_{3}},~p_1,~p_2, \cdots,~p_{t},~p_{t+1},\cdots,~p_{t+1}}$
such that the number of its parts consisting of $p_{t+1}$ vertices is $n+1-t$.
It is clear that  $G_2$ is a subgraph of $G$. Then $\theta'(G_2)\leq \theta'(G)= n-t$.
We can partition the vertex set $V(G_2)$ into $n-t$ (or fewer) subsets such that each subset induces a planar subgraph.
The total number of vertices of $G_2$ is
$$
\sum^{t+1} _{i=0} p_i +(n+1-t-1)p_{t+1}
> 2(n-t)+(n-t)p_{t+1}
=(n-t)(2+p_{t+1}).
$$
This implies that in any partition of $V(G_2)$ into $n-t$ (or fewer) subsets, there is a
 subset which contains at least  $3+p_{t+1}$ vertices.
Since $3+p_{t+1}\ge 3+ p_1\ge 7$, the subset, say $V'$, contains $m\ge 7$ vertices. Since any part in $G_2$ has at most $p_{t+1}$ vertices,  the subgraph induced by $V'$ is not $K_{1, m-1}$ ,  $K_{2, m-2}$, or $K_{1,1, m-2}$.
From Lemma \ref{lem1}, the subgraph is a non-planar graph, a contradiction.

\medskip
\textit{ Case 2: Suppose $\sum^{t+1} _{i=0} p_i \le 2(n+1-(t+1))=2(n-t)$. }
\medskip

It follows from the second inequation in $(*)$ that
$$\sum^{t+2} _{i=0} p_i =p_{t+2}+ \sum^{t+1} _{i=0} p_i > 2(n-t-1)=2(n+1-(t+2)),$$
 which implies
$$ n+1-max\{ j| \sum^{j} _{i=0} p_i \leq 2(n+1-j) \}=n+1-(t+1)=n-t=\theta'(G_1).$$
Thus, we need to prove  $\theta'(G)=\theta'(G_1)$ in this case.

Let $V'$ be a subset of $V(G)$ consisting of all vertices of the first $(k_1+k_2+k_3+t+1)$ parts in $G$.
 Since $\sum^{t+1} _{i=0} p_i \le 2(n-t)$, we can exhaust the set $V'$ by
 adding at most two vertices of  $V'$  to each of the sets
$V_{t+2}, V_{t+3}, \cdots, V_{n+1}$. Every one of the resulting subsets induces a planar graph isomorphic to
$K_{2, p_s}$,  $K_{1, 1, p_s}$, $K_{1, p_s}$ or an  $I_{p_s}$-set,  where
$t+2 \leq s \leq n+1$.
Thus, $\theta'(G)\leq n+1-(t+1)=n-t$.  Since $G_1$ is a subgraph of $G$, we have $\theta'(G)\ge \theta'(G_1)=n-t$. Therefore,
$$~~~~~~~~~~~~~~~~~~~~~~~~~~~~\theta'(G)=\theta'(G_1)=n-t=n+1-max\{ j| \sum^{j} _{i=0} p_i \leq 2(n+1-j) \}. ~~~~~~~~~~~~~~~~~~~~~~~~~~~~~~~~~~~~~\qed $$

\begin{lem}\label{lem3}
Let $G, p_0$ and $n$ be defined as above. If $p_0 > 2n$, then

$$\theta'(G)=\left\{
                \begin{array}{ll}
                  n+N(k_1+k_3-2n, k_2+k_3), & \hbox{if $k_1+k_3 \geq 2n$;} \\
                  n+N(\varepsilon, \frac{1}{2}(p_0-2n-\varepsilon)), & \hbox{otherwise,}
                  \end{array}
              \right. $$
where  $k_1+k_3\equiv \varepsilon $ (mod 2).

\end{lem}

\pf  Let $G_1$ be a complete multipartite graph $K_{1^{k_{1}},~2^{k_{2}},~3^{k_{3}},~4,~4, \cdots, 4}$ such that the number of parts consisting of four vertices is $n$.
Then $G_1$ is a subgraph of $G$ and hence $\theta'(G_1)\le \theta'(G)$.
There are $k_2+k_3+2n$ $I_2$-sets in $G_1$,
since each part containing two or three vertices contributes one $I_2$-set and each part with four vertices can  contribute two $I_2$-sets.

For the complete multipartite graph $G=K_{1^{k_{1}},~2^{k_{2}},~3^{k_{3}},~p_1,~p_2, \cdots,~p_n}$,
there are $k_3$  parts containing exact three independent vertices, i.e., an $I_3$-set.
We divide each of the $k_3$  $I_3$-sets into an  $I_1$-set and an $I_2$-set.
Let  $\mathcal{I}_1$ be such a set in $G$ that consists of $k_3$ $I_1$-sets obtained  by dividing $k_3$ $I_3$-sets  and
$k_1$ $I_1$-sets which are the parts with  only one vertex, and let
 $\mathcal{I}_2$ be such a set in $G$ that consists of $k_3$ $I_2$-sets obtained  by dividing $k_3$ $I_3$-sets and
$k_2$ $I_2$-sets which are the parts with two vertices.

According to the definition of $G$,  the sets $V_1$, $V_2$, $\cdots$, $V_n$ are the last $n$ subsets of the vertex partition of $G$ satisfying
$|V_i|=p_i$ ($1\le i\le n $).
Since $p_0 > 2n$, we can add two elements of $\mathcal{I}_1$ or one of $\mathcal{I}_2$ into one of 
$\{V_i| 1\le i\le n \}$ each time 
until  all $V_i$ are exhausted. 
Our proof consists of the following two parts.

\medskip
{\bf\textit{  Part I:  $k_1+k_3 \ge 2n$, i.e., $\mathcal{I}_1$ can exhaust all  $V_i$ $(1\le i\le n)$.}}
\medskip

\textit{ Case 1.1.   $k_1+k_3-2n \ge 3(k_2+k_3)$ and further $k_1-3k_2-2k_3-2n \equiv 0$ (mod 4).}
\medskip

Let $s=\frac{1}{4}(k_1-3k_2-2k_3-2n)$ and $t=k_2+k_3+s$. First, we add two elements of $\mathcal{I}_1$ to one  of $\{V_i| 1\le i\le n \}$ 
  each time until all $V_i$ are exhausted, resulting in $n$ sets. Each of 
the $n$ sets induces a planar graph $K_{1,1, p_i}$ ($1\le i\le n $). Then, 
since $k_1+k_3-2n \ge 3(k_2+k_3)$,
we can use up all elements of
 $\mathcal{I}_2$ by repeatedly adding three of the remaining elements of  $\mathcal{I}_1$ to one element of $\mathcal{I}_2$, and obtain $k_2+k_3 $ sets. 
 Every one of the $k_2+k_3 $ sets
 induces a planar graph  $K_{1, 1, 1, 2}$ or $K_{1, 1,  3}$.  Finally, from $k_1-3k_2-2k_3-2n \equiv 0$ (mod 4), the unused elements of $\mathcal{I}_1$ can be divided
into $s$ sets each of which consists of four vertices and induces a planar graph $K_4$.
Thus,  $V(G)$ is partitioned into $n+t$ subsets such that each subset induces a planar graph. Therefore $\theta'(G)\le n+t$.
To prove $\theta'(G)= n+t$, it suffices to prove $\theta'(G)\ge n+t$.

\medskip
\textit{Claim 1:}~$\theta'(G)\ge n+t$.
\medskip

\pf  Since $G_1$ is a subgraph of $G$,
$\theta'(G)\ge \theta'(G_1)$.  It is sufficient to prove  $\theta'(G_1)\ge n+t$.
(In fact, we can  prove $\theta'(G_1)= n+t$.)

Suppose by contradiction $\theta'(G_1)\le n+t-1$. Then we can partition $V(G_1)$ into at most $n+t-1$ subsets so that each subset
induces a planar graph.     Denote by $\mathcal{A}$ the set consisting of all partitioned subsets of $V(G_1)$.
Clearly,
$$k_1+2k_2+3k_3+4n=4(n+t-1)+(k_2+k_3+2n+4).$$

If there is an element of $\mathcal{A}$ containing $m\ge 7$ vertices, then since each part of $G_1$ has at most four vertices,
 the subgraph of $G_1$ induced by the element  with $m$ vertices
is not  $K_{1, m-1} $,  $K_{2, m-2}$, or $K_{1, 1,  m-2}$,  and hence is a non-planar graph by Lemma \ref{lem1}, a contradiction. Thus  every element of $\mathcal{A}$
has at most six vertices.

We assume that  $\mathcal{A}$ contains $k\ge 0$ elements with six vertices.  Then there are at least
$(k_2+k_3+2n+4-2k)$ elements with five vertices in $\mathcal{A}$: If not, then the number of elements in $\mathcal{A}$ with
five vertices is at most $k_2+k_3+2n+3-2k$. 
The total number of vertices of $G_1$ is at most
\begin{align*}
&~~~6k+5(k_2+k_3+2n+3-2k)
+ 4(n+t-1-k-(k_2+k_3+2n+3-2k))\\
&=k_1+2k_2+3k_3+4n-1,
\end{align*}
which contradicts $|V(G_1)|=k_1+2k_2+3k_3+4n$.

From the definition of  $\mathcal{A}$, the subgraph in $G_1$ induced by any element of $\mathcal{A}$ is a planar graph. Since all parts of $G_1$ have at most four vertices,  the subgraph
induced by an element of $\mathcal{A}$ with six vertices is $K_{2, 4}$, $K_{1, 1, 4}$ or $K_{2, 2, 2}$, each having at least two $I_2$-sets, and
the subgraph
induced by an element of $\mathcal{A}$ with five vertices is $K_{1,1,1,2}$, $K_{1,2,2}$, $K_{1,1,3}$, $K_{2,3}$, or $K_{1,4}$, each having at
least one $I_2$-set.
Then there are at least
$$2k+(k_2+k_3+2n+4-2k)
=k_2+k_3+2n+4$$
$I_2$-sets in $G_1$,  a  contradiction with the fact that
$G_1$ contains $k_2+k_3+2n$ $I_2$-sets.
Thus  $\theta'(G_1)\ge n+t$ and so $\theta'(G)\ge \theta'(G_1) \ge n+t$.
\qed

\bigskip
Since $k_1-3k_2-2k_3-2n \equiv 0$ (mod 4), we have
\begin{align*}
 \theta'(G)&= n+t=n+k_2+k_3+\frac{1}{4}(k_1+k_3-2n-3(k_2+k_3))\\
&=n+N(k_1+k_3-2n, k_2+k_3).
\end{align*}

\medskip
\textit{ Case 1.2.   $k_1+k_3-2n \ge 3(k_2+k_3)$ and further $k_1-3k_2-2k_3-2n \equiv \gamma$ (mod 4) with $\gamma$=1, 2, or 3.}
\medskip

Let $s=\frac{1}{4}(k_1-3k_2-2k_3-2n -\gamma)$ and $t=k_2+k_3+s+1$. 
We partition $V(G)$ into $n+t$ subsets below.
First,  we  add two elements of $\mathcal{I}_1$ to one of $\{V_i| 1\le i\le n\}$
each time until all  $V_i$ are exhausted.
 The $n$ resulting subsets induce $n$ planar graphs each isomorphic to $K_{1,1,p_i}$ ($1\le i\le n$).
 Then, from  $k_1+k_3-2n \ge 3(k_2+k_3)$,
we add three of the remaining elements of $\mathcal{I}_1$ to one element of $\mathcal{I}_2$ each time until all elements of
$\mathcal{I}_2$ are exhausted, resulting in $k_2+k_3$ subsets. The graphs induced by the $k_2+k_3$  subsets
are either $K_{1,1,1,2}$ or $K_{1,1, 3}$.
Finally, since $k_1-3k_2-2k_3-2n \equiv \gamma$ (mod 4) with $\gamma$=1, 2, or 3, the unused elements of $\mathcal{I}_1$ can be divided
into $s+1$ subsets such that one subset consists of $\gamma$ vertices and $s$ subsets all consist of four vertices.
All graphs induced by the $s+1$ subsets are planar graphs since  their orders are at most 4. 
Thus, $V(G)$ is divided into $n+t$ subsets such that each subset induces a planar graph,  which implies  $\theta'(G)\le n+t$.
Similar to Claim 1, we can prove $\theta'(G)\ge n+t$. Therefore $\theta'(G)= n+t$. 

It follows from $k_1-3k_2-2k_3-2n \equiv \gamma$ (mod 4) with $\gamma$=1, 2, or 3 that 
\begin{align*}
\theta'(G)&=n+t=n+k_2+k_3+\frac{1}{4}(k_1-3k_2-2k_3-2n-\gamma)+1\\
&=n+k_2+k_3+\frac{1}{4}(k_1+k_3-2n-3(k_2+k_3))+\frac{1}{4}(4-\gamma)\\
&=n+k_2+k_3+\frac{1}{4}\left\lceil k_1+k_3-2n-3(k_2+k_3)\right\rceil\\
&=n+N(k_1+k_3-2n,k_2+k_3).
\end{align*}

\medskip
\textit{ Case 1.3. $k_1+k_3-2n<3(k_2+k_3)$ with $k_1+k_3-2n\equiv \gamma_1$ (mod 3) and $k_2+k_3-\frac{1}{3}(k_1+k_3-2n-\gamma_1)\equiv \gamma_2$ (mod 3), where
$0 <\gamma_1+\gamma_2<4$.}
\medskip

Let $s=\frac{1}{3}(k_1+k_3-2n-\gamma_1)$ and $t=\frac{1}{3}(k_2+k_3-s-\gamma_2)+s+1$. We first add two elements of $\mathcal{I}_1$ into one of the sets $V_i$ ($1\le i\le n$)
each time until all the sets $V_i$ are exhausted, resulting in $n$ subsets.
The resulting subsets induce $n$ planar graphs 
each isomorphic to $K_{1,1,p_i}$($1\le i\le n$). Then, since $k_1+k_3-2n<3(k_2+k_3)$ and $k_1+k_3-2n\equiv \gamma_1$ (mod 3),
we  add three of the remaining elements of $\mathcal{I}_1$ into one element of $\mathcal{I}_2$ each time until \ $\mathcal{I}_1$ only remains $\gamma_1$ elements, and obtain $s=\frac{1}{3}(k_1+k_3-2n-\gamma_1)$ subsets each with five vertices.
 Every one of the $s$ subsets
induces a planar graph $K_{1,1,1,2}$ or $K_{1,1,3}$. Next, from $k_2+k_3-s\equiv \gamma_2$ (mod 3), the remaining elements of $\mathcal{I}_2$ can be divided
into $\frac{1}{3}(k_2+k_3-s-\gamma_2)$ subsets each of which consists of six vertices,
and has exact $\gamma_2$ elements left. Each of the $\frac{1}{3}(k_2+k_3-s-\gamma_2)$ subsets
 induces a planar graph $K_{2,2,2}$. 
Finally, the remaining $\gamma_1$ elements of $\mathcal{I}_1$ and  $\gamma_2$ elements of $\mathcal{I}_2$ together form a subset, which can
 induce a planar graph since $0 <\gamma_1+\gamma_2<4$.
Thus, we partition $V(G)$ into $n+t$ subsets so that each subset induces a planar graph. Therefore $\theta'(G)\le n+t$.
We need to prove $\theta'(G)\ge n+t$ in order to prove $\theta'(G)= n+t$.

\medskip
\textit{Claim 2:}~$\theta'(G)\ge n+t$.
\medskip

\pf  Since $G_1$ is a subgraph of $G$, we deduce  
$\theta'(G)\ge \theta'(G_1)$. It is sufficient to prove $\theta'(G_1)\ge n+t$.
Suppose by contradiction that $\theta'(G_1)\le n+t-1$. Then $V(G_1)$ can be partitioned into at most $n+t-1$ subsets such that each subset induces a planar graph. 
We use $\mathcal{A}$ to denote the set  of all partitioned subsets.
It is clear that 
$$k_1+2k_2+3k_3+4n=5(n+t-1)+\frac{1}{9}(3k_2+2k_3-k_1+2n)+n+\frac{10}{9}\gamma_1+\frac{5}{3}\gamma_2.$$

If $\mathcal{A}$ contains an element with $m\ge 7$ vertices, then since all parts of $G_1$ have at most four vertices,
the subgraph in $G_1$ induced by the element with $m$
is not $K_{1,m-1}$, $K_{2,m-2}$ or $K_{1,1,m-2}$, and hence is a non-planar graph by Lemma \ref{lem1}, a contradiction. Thus, every element of $\mathcal{A}$
has at most six vertices.

Let  $m=\frac{1}{9}(3k_2+2k_3-k_1+2n)+n+\frac{10}{9}\gamma_1+\frac{5}{3}\gamma_2$. Then 
there are at least $m$ elements in $\mathcal{A}$ which all contain six vertices: If not, the number of elements in $\mathcal{A}$ containing six vertices is at most $(m-1).$
Since other elements of $\mathcal{A}$ have
at most five vertices,  the total number of vertices of $G_1$ is not more than
$$6(m-1)+5(n+t-1-(m-1))=k_1+2k_2+3k_3+4n-1.$$
which contradicts $|V(G_1)|=k_1+2k_2+3k_3+4n$.

Therefore we assume $\mathcal{A}$
contains $m+k$ ($k\ge 0$) elements with six vertices.  Then the number of elements in $\mathcal{A}$ with five vertices is at least
$(\frac{1}{3}k_1+\frac{1}{3}k_3-\frac{2}{3}n-\frac{4}{3}\gamma_1-2\gamma_2-2k)$: If not, the number is at most $(\frac{1}{3}k_1+\frac{1}{3}k_3-\frac{2}{3}n-\frac{4}{3}\gamma_1-2\gamma_2-2k-1)$. The other elements in $\mathcal{A}$ different from the elements with six and five vertices
have at most four vertices, and 
$$n+t-1-(m+k)-(\frac{1}{3}k_1+\frac{1}{3}k_3-\frac{2}{3}n-\frac{4}{3}\gamma_1-2\gamma_2-2k-1)=k+1.$$
Then, the total number of vertices of $G_1$ is not  greater than
\begin{align*}
&6(m+k)+5(\frac{1}{3}k_1+\frac{1}{3}k_3-\frac{2}{3}n-\frac{4}{3}\gamma_1-2\gamma_2-2k-1)+4(k+1)\\
=&k_1+2k_2+3k_3+4n-1,
\end{align*}
which is a contradiction with $|V(G_1)|$.

As  all parts of $G_1$ have at most four vertices,  the subgraph
induced by the element of $\mathcal{A}$ with six vertices is $K_{2, 4}$, $K_{1, 1, 4}$ or $K_{2, 2, 2}$.  Clearly,  both $K_{2, 4}$ and $K_{2, 2, 2}$  contain three $I_2$-sets, while $K_{1, 1, 4}$ contains two $I_2$-sets. Since the number of the parts in $G_1$ with four vertices is $n$, the number of  the induced subgraph $K_{1, 1, 4}$ is at most $n$. It is easily checked that
 the subgraph
induced by the element of $\mathcal{A}$ with five vertices is $K_{1,1,1,2}$, $K_{1,2,2}$, $K_{1,1,3}$, $K_{2,3}$, or $K_{1,4}$, each having at
least one $I_2$-set.
Then, in $G_1$ there are at least
\begin{align*}
   &2n+3(m+k-n)+(\frac{1}{3}k_1+\frac{1}{3}k_3-\frac{2}{3}n-\frac{4}{3}\gamma_1-2\gamma_2-2k)\\
  =&k_2+k_3+2n+2\gamma_1+3\gamma_2+k
\end{align*}
$I_2$-sets, which contradicts the fact that
$G_1$ contains $k_2+k_3+2n$ $I_2$-sets since $2\gamma_1+3\gamma_2+k>0$.
So $\theta'(G_1)\ge n+t$ and hence $\theta'(G)\ge \theta'(G_1)\ge n+t$.
\qed

\bigskip
From $k_1+k_3-2n\equiv \gamma_1$ (mod 3) and $k_2+k_3-\frac{1}{3}(k_1+k_3-2n-\gamma_1)\equiv \gamma_2$ (mod 3), we obtain
\begin{align*}
\theta'(G)&=n+t=n+\frac{2}{9}k_1+\frac{1}{3}k_2+\frac{5}{9}k_3-\frac{4}{9}n-\frac{2}{9}\gamma_1-\frac{1}{3}\gamma_2+1\\
&=n+\frac{1}{3}(k_1+k_3-2n-\gamma_1)+\frac{1}{3}(k_2+k_3-\frac{1}{3}(k_1+k_3-2n-\gamma_1)-\gamma_2)+1\\
&=n+\lfloor \frac{1}{3}(k_1+k_3-2n)\rfloor+\lfloor\frac{1}{3}(k_2+k_3-\lfloor \frac{1}{3}(k_1+k_3-2n)\rfloor)\rfloor+1\\
&=n+N(k_1+k_3-2n,k_2+k_3)
\end{align*}

\medskip
\textit{ Case 1.4. $k_1+k_3-2n<3(k_2+k_3)$ with $k_1+k_3-2n\equiv 2$ (mod 3) and $k_2+k_3-\frac{1}{3}(k_1+k_3-2n-2)\equiv 2$ (mod 3).}
\medskip

Let $s=\frac{1}{3}(k_1+k_3-2n-2)$ and  $t=\frac{1}{3}(k_2+k_3-s-2)+s+2$. First,  we  add two elements of $\mathcal{I}_1$ into one element of $\{V_i| 1\le i\le n\}$
each time until all the sets $V_i$ are exhausted, resulting in  $n$ subset.
Each of the resulting subsets induces a planar graph $K_{1,1,p_i}$($1\le i\le n$). Second, since  $k_1+k_3-2n<3(k_2+k_3)$ and $k_1+k_3-2n\equiv 2$ (mod 3),
we can repeatedly add three of the remaining elements of $\mathcal{I}_1$ into one element of $\mathcal{I}_2$  until  $\mathcal{I}_1$ remains exact  $2$ elements, and obtain $s$ subsets with five vertices.  Each of the $s$  subsets
induces a planar graph $K_{1,1,1,2}$ or $K_{1,1,3}$. Third, from $k_2+k_3-s\equiv 2$ (mod 3), the remaining elements of $\mathcal{I}_2$ can be divided to $\frac{1}{3}(k_2+k_3-s-2)$ subsets of six vertices, and then have exact $2$ elements left. 
Every one of the $\frac{1}{3}(k_2+k_3-s-2)$ subsets
 induces a planar graph $K_{2,2,2}$.  Finally,
the remaining elements of $\mathcal{I}_1$ and elements of $\mathcal{I}_2$
sum up to six vertices and induce a non-planar subgraph $K_{1,1,2,2}$, $K_{1,2,3}$ or $K_{3,3}$
So we  partition these six vertices into any two empty subsets, so that each subset can induce a planar graph. 
Thus $V(G)$ is divided into $n+t$ subsets such that each subset induces a planar graph. This implies $\theta'(G)\le n+t$. Similar to Claim 2, we can prove $\theta'(G)\ge n+t$.
Therefore $\theta'(G)=n+t$. 

Since $k_1+k_3-2n\equiv 2$ (mod 3) and $k_2+k_3-\frac{1}{3}(k_1+k_3-2n-2)\equiv 2$ (mod 3), 
\begin{align*}
\theta'(G)&=n+t=n+\frac{2}{9}k_1+\frac{1}{3}k_2+\frac{5}{9}k_3-\frac{4}{9}n+\frac{8}{9}\\
&=n+\frac{1}{3}(k_1+k_3-2n-2)+\frac{1}{3}(k_1+k_3-\frac{1}{3}(k_1+k_3-2n-2)-2)+2\\
&=n+\lfloor \frac{1}{3}(k_1+k_3-2n) \rfloor+\lfloor\frac{1}{3}(k_2+k_3-\lfloor\frac{1}{3}(k_1+k_3-2n)\rfloor)\rfloor+2\\
&=n+N(k_1+k_3-2n,k_2+k_3)
\end{align*}

\medskip
\textit{ Case 1.5. $k_1+k_3-2n<3(k_2+k_3)$ with $k_1+k_3-2n\equiv 0$ (mod 3) and  $k_2+k_3-\frac{1}{3}(k_1+k_3-2n)\equiv 0$ (mod 3).}
\medskip

Similar to Case 1.3, we can prove $\theta'(G)=n+N(k_1+k_3-2n,k_2+k_3)$. As it is easier than Case 1.3, we omit all details of the proof.

\bigskip
{\bf\textit{  Part II:  $k_1+k_3 < 2n$, i.e., $\mathcal{I}_1$ cannot exhaust all  $V_i$  $(1\le i\le n)$.}}
\medskip

\medskip
\textit{ Case 2.1. $k_1+k_3\equiv \varepsilon$ (mod 2) and $\frac{1}{2}(p_0-2n-\varepsilon)\equiv \gamma$ (mod 3),  where $\varepsilon +\gamma >0$.}
\medskip

Let $t=\frac{1}{3}(\frac{1}{2}(p_0-2n-\varepsilon)-\gamma)+1$. We first add two elements of $\mathcal{I}_1$ into one element of the set $\{V_i| 1\le i \le n \}$
each time until $\mathcal{I}_1$ remains exact $\varepsilon$  elements. 
The set $\{V_i| 1\le i \le n \}$
has $\frac{1}{2}(k_1+k_3-\varepsilon)$ elements used 
and  $n-\frac{1}{2}(k_1+k_3-\varepsilon)$ elements left. 
Then, we add one element of $\mathcal{I}_2$ into one of the remaining elements of $\{V_i| 1\le i \le n \}$
each time until  the set $\{V_i| 1\le i \le n \}$ are exhausted. We obtain $n$ subsets in which each set
induces a planar graph $K_{1,1,p_i}$ or $K_{2,p_i}$ ($1\le i\le n$).
Here
$\mathcal{I}_2$ remains $\frac{1}{2}(p_0-2n-\varepsilon)$ elements.
From $\frac{1}{2}(p_0-2n-\varepsilon)\equiv \gamma$ (mod 3), the remaining elements of $\mathcal{I}_2$ can be divided
into $\frac{1}{3}(\frac{1}{2}(p_0-2n-\varepsilon)-\gamma)$ subsets each of which consists of six vertices and induces a planar graph $K_{2,2,2}$, and then have exact $\gamma$ elements left.
Finally, the  $\varepsilon$ remaining elements of $\mathcal{I}_1$ together with the $\gamma$ remaining elements of $\mathcal{I}_2$ can induce a planar graph.
Thus, we partition $V(G)$ into $n+t$ subsets such that each subset induces a planar graph. Therefore $\theta'(G)\le n+t$.
To prove $\theta'(G)= n+t$, it suffices to prove $\theta'(G)\ge n+t$.

\medskip
\textit{Claim 3:}~$\theta'(G)\ge n+t$.
\medskip

\pf Since $\theta'(G)\ge \theta'(G_1)$, we only need  to prove $\theta'(G_1)\ge n+t$.
Suppose by contradiction $\theta'(G_1)\le n+t-1$. 
Then $V(G_1)$ can be partitioned into at most $n+t-1$ subsets such that every subset
induces a planar graph. Denote by $\mathcal{A}$ the set consisting of all partitioned subsets.
It can be computed that 
$$k_1+2k_2+3k_3+4n=6(n+t-1)+\varepsilon+2\gamma.$$
Since $\varepsilon+\gamma >0$ and hence $\varepsilon+2\gamma >0$, 
there is at least one  element of $\mathcal{A}$ containing $m\ge 7$ vertices. The subgraph in $G_1$ induced by the element with $m$ vertices 
is not $K_{1,m-1}$, $K_{2,m-2}$ or $K_{1,1,m-2}$, and hence is a non-planar graph according to Lemma \ref{lem1}, a contradiction.
Thus, $\theta'(G)\ge \theta'(G_1)\ge n+t$.
\qed

\medskip
From $k_1+k_3\equiv \varepsilon$ (mod 2) and $\frac{1}{2}(p_0-2n-\varepsilon)\equiv \gamma$ (mod 3),
we deduce
\begin{align*}
\theta'(G)&=n+t=n+\frac{1}{3}(\frac{1}{2}(p_0-2n-\varepsilon)-\gamma)+1\\
&=n+\lfloor\ \frac{\varepsilon}{3} \rfloor+\lfloor \frac{1}{3}(\frac{1}{2}(p_0-2n-\varepsilon)-\lfloor \frac{\varepsilon}{3} \rfloor)+1\\
&=n+N(\varepsilon,\frac{1}{2}(p_0-2n-\varepsilon)) \\
\end{align*}

\textit{ Case 2.2. $k_1+k_3\equiv 0$ (mod 2) and $\frac{1}{2}(p_0-2n)\equiv 0$ (mod 3).}
\medskip

Let $t=\frac{1}{6}(p_0-2n)$. Similar to Case 2.1, we can prove $\theta'(G)= n+t$.
Thus, 
\begin{align*}
~~~~~~~~~~~~~~~~~~~~~~~~~~~~
\theta'(G)&=n+t=n+\frac{1}{6}(p_0-2n)\\
&=n+\lfloor \frac{0}{3}\rfloor+\lfloor\frac{1}{3}(\frac{1}{2}(p_0-2n)-\lfloor \frac{0}{3}\rfloor)\rfloor\\
&=n+N(0,\frac{1}{2}(p_0-2n))~~~~~~~~~~~~~~~~~
~~~~~~~~~~~~~~~~~~~~~~~~~~~~~~~~~~~~~\qed
\end{align*}

\bigskip

Now we are ready to prove our main theorem.

{\bf Proof of  Theorem 1.1}. Theorem 1.1 $(a)$  follows from Lemma \ref{lem2} and
Theorem 1.1 $(b)$  from Lemma \ref{lem3}.

\medskip
From Theorem 1.1, we easily deduce the following corollary. 
\medskip

{\bf Corollary 1.2.} For a  complete graph $K_n$, then
$\theta'(K_n)=\lceil   \frac{n}{4}  \rceil.$

\bigskip

{\bf Remark.} 
From the proof of  Lemma \ref{lem3}, we know that the number $N(k_1,k_2)$ is, in fact, the point-thickness $\theta'(K_{1^{k_{1}}, 2^{k_{2}}})$ of 
the complete multipartite graph $K_{1^{k_{1}}, 2^{k_{2}}}$.
For a complete multipartite graph $G=K_{1^{k_{1}}, 2^{k_{2}}, 3^{k_{3}}, p_1, p_2, \cdots, p_n}$
with $4 \leq p_1\leq  p_2\leq \cdots \leq p_n$, when $p_0=k_1+2k_2+3k_3>2n$, its point-thickness 
$\theta'(G)$  depends  on $n$ and is independent of the specific values of all $p_i$ $(1\le i\le n)$.

\bigskip

{\bf Declaration of competing interest.}
The authors declare that they have no known competing financial interests or personal relationships that could have
appeared to influence the work reported in this paper.

\bigskip

{\bf Acknowledgements.}  Wenzhong Liu was
partially supported by the National Natural Science Foundation of China (No.12271251) and
the Fundamental Research Funds for the Central Universities (No.NC2024007).

\end{document}